\documentclass[12pt]{article}
\usepackage[margin=1.5in]{geometry} 
\usepackage{amsmath,amsthm,amssymb,amsfonts}
\usepackage{changepage}
\usepackage{tikz}
\usepackage[shortlabels]{enumitem}
\usetikzlibrary{matrix}

\newcommand{\Z}{\mathbb{Z}}

\newcommand{\inv}{^{-1}}

\newcommand{\res}{\upharpoonright}

\newtheorem{theorem}{Theorem}

\newtheorem{lemma}{Lemma} 
  
\newtheorem{prob}{Problem} 
 
\begin{document}
 
\title{Borel Vizing's Theorem for 2-Ended Groups}
\author{Felix Weilacher}
\maketitle

\begin{abstract}
\noindent We show that Vizing's Theorem holds in the Borel context for graphs induced by actions of 2-ended groups, and ask whether it holds more generally for everywhere two ended Borel graphs.
\end{abstract}


\section{Introduction}\label{sec:intro}


For a graph $G$ on a set $X$, let $\chi'(G)$ denote the edge chromatic number of $G$. That is, the smallest cardinal $k$ such that there exists a function assigning each edge in $G$ to an element of $k$ such that any two edges incident on the same vertex are assigned different elements. Such a function is called a $k$-edge coloring of $G$, and the elements of $k$ are called colors. If $G$ is a Borel graph on a standard Borel space $X$, let $\chi_B'(G)$ denote the Borel edge chromatic number of $G$. That is, the smallest cardinal $k$ as above, but where only Borel (as functions) colorings are allowed.

A classical theorem of Vizing states that if $G$ is a graph of maximum degree $d \in \omega$, $\chi'(G) \leq d+1$. Note the tirival lower bound $\chi'(G) \geq d$, so that Vizing's theorem implies $\chi'(G) \in \{d,d+1\}$. We are interested in generalizations of this theorem to the Borel context. Marks has shown \cite{M15} that the direct generalization fails, but on the other hand, Grebik and Pikhurko have shown \cite{GP} that the generalization holds if `Borel' is weakened to `$\mu$-measurable' for some Borel probability measure $\mu$ on $X$ which is $G$-invariant.

Recently, Weilacher has shown \cite{W20} that some combinatorial bounds which hold in the measurable context but not generally in the Borel context can still be salvaged in the Borel context with an additional assumption: that every connected component of $G$ has two ends. In the spirit of this, we show in this note the following:
\begin{theorem}\label{th:main}
Let $\Gamma$ be a marked group with two ends, say with $d$ generators. Let $G$ be the shift graph of $\Gamma$, so that $G$ is $d$ regular. Then $\chi_B'(G) \leq d+1$.
\end{theorem}

One purpose of this note is to pose the question of whether the assumption in Theorem \ref{th:main} that $G$ be generated by a group action is necessary, which seems to be open.

\begin{prob}
Let $G$ be a Borel graph of maximum degree $d$ such that every connected component of $G$ has two ends. Is $\chi'_B(G) \leq d+1$?
\end{prob}

\section{Proof}

In this section we prove Theorem \ref{th:main}. The proof is very simple and intuitive, but it takes some time to write down all of the details

\begin{proof}

Fix a two ended marked group $\Gamma$ with symmetric generating set $S$ of size $d$. It is well known that since $\Gamma$ has two ends, there is a finite normal subgroup $\Delta \leq \Gamma$ such that $\Gamma/\Delta \cong \Z$ or $D_\infty = \langle a,b \mid a^2=b^2=\textnormal{id} \rangle$. Let us start with the former case for ease of notation. The latter case can be handled in the same way and will be addressed at the end of this section.

Partition $S$ as $S = \bigsqcup_{n \in \Z} S_n$, where $S_n = \{\gamma \in S \mid \overline{\gamma} = n\}$, where $\overline{\gamma}$ denotes the image of $\gamma$ in the quotient $\Gamma/\Delta$ identified with $\Z$. Note that $S_{-n} = S_n\inv = \{\gamma\inv \mid \gamma \in S_n\}$ for each $n$. Let $d_n = |S_n|$ for each $n$, so that $\sum_n d_n = d$.

Let $G$ be the shift graph of $\Gamma$ with vertex set $X$. Let $Y$ be the standard Borel space of $\Delta$-orbits of $X$. The action of $\Gamma$ on $X$ descends to an action of $\Z$ on $Y$. Let $H$ denote the Borel multigraph on $Y$ defined by placing an edge between $y$ and $n \cdot y$ for each $\gamma \in S_n$ for $n > 0$ and $y \in Y$. In other words, the number of edges between $y$ and $n \cdot y$ in $H$ is always $d_n$ for $n \neq 0$. Let $k = \sum_{n > 0} d_n$, so that $H$ is $2k$-regular. We claim that $H$ admits a Borel $2k+1$-edge coloring. Since $H$ is generated by an action of $\Z$, it suffices to just prove our main Theorem when $\Gamma = \Z$ and multiplicity for the generators is allowed:

\begin{lemma}\label{lem:Z}
Theorem \ref{th:main} holds when $\Gamma = \Z$, even when multiplicity for generators is allowed. 
\end{lemma}

It should be noted that Vizing's Theorem fails in general for multigraphs, (although there is a generalization which still holds) so this lemma is somewhat surprising. We now prove it:

\begin{proof}
Keeping in line with the notation established so far, let $H$ be the shift graph of $\Z$ with our generating set with vertex set $Y$. Let $0 < n_1 \leq n_2 \leq \cdots \leq n_k$ list the positive generators with multiplicity, so that $H$ is $2k$-regular.

Let $G'$ be the graph on $Y$ induced by the action of $\Z$ with usual generators $\pm 1$.

We first show the following: Let $A \subset Y$ be a Borel subset of our shift graph such that the induced subgraph $G' \res A$ has connected components all of size at least $2$, and such that $A$ is \textit{recurrent}. That is, for each $x \in Y$, there are $m,l > 0$ such that $m \cdot x \in A$ and $(-l) \cdot x \in A$. Then there is a Borel 3-edge coloring of $G'$, say using the colors 1,2, and 3, such that the color 3 only appears on edges between points in $A$.

First, by the recurrence of $A$, we can clearly find a Borel $G'$-independent recurrent set $B \subset Y$ such that if $x,y \in B$ are distinct points in the same connected component of $G'$, the unique path between them in $G'$ passes through $A$. Suppose $x,y \in B$ such that $y = N \cdot x$ for some $N > 0$ and there are no points of $B$ between $x$ and $y$ in the graph $G'$. 

We need to color the edges between $x$ and $y$. If $N$ is even, we color the edge $(m \cdot x, (m+1) \cdot x)$ with the color 1 for $m$ even and the color 2 for $m$ odd for $0 \leq m < N$. If $N$ is odd, let $0 < M < N-1$ be minimal with the property that $M \cdot x \in A$. This exists by definition of $B$, and then $(M+1) \cdot x \in A$ by definition of $A$. Accordingly, we color the edge $(m \cdot x,(m+1) \cdot x)$ with the color 1 for $m$ even and the color 2 for $m$ odd for $0 \leq m < M$, then color the edge $(M \cdot x, (M+1 \cdot x))$ with the color 3, then color the edge $(m \cdot x,(m+1) \cdot x)$ with the color $1$ for $m$ odd and the color 2 for $m$ even for $M < m < N$. Note that for each $x \in B$, $(x, 1 \cdot x)$ has color 1 and $(-1 \cdot x, x)$ has color 2, so we do indeed end up with a coloring. Furthermore, the color 3 was clearly only ever used for edges between points of $A$.

Now, returning to our original goal, we begin by partitioning $Y$ into $k$ Borel recurrent sets $A_1,\ldots,A_k$ such that for each $i$, the connected components of $G' \res A_i$ each have size at least $2n_k$. By a result from \cite{KST}, we can start by finding a Borel maximal $2n_k$-discrete set $B$. In particular $B$ will be recurrent. We can then partition $B$ into $k$ many recurrent sets $B_1,\ldots,B_k$ using, for example, that same result. Now, for each $i$ and each $x \in B_i$, there will be a smallest $N > 0$ such that $N \cdot x \in B$. We then include $m \cdot x \in A_i$ for each $0 \leq m < N$. This clearly works.

Now, fix one of our generators $n_i$. Consider only the edges in $H$ corresponding to this generator, and call the resulting (simple) graph $H_i$. Abusing language slightly, observe that since $n_i \leq n_k$, our set $A_i$ will be recurrent for the graph $H_i$ (more precisely, for each $y \in Y$, there are $m,l > 0$ such that $mn_i \cdot y$ and $(-l)n_i \cdot y \in A_i$), and all the connected components of $H_i \res A_i$ will have size at least 2. It follows from the statement of paragraph 3 of this proof that we can Borel edge color $H_i$, say using the colors $2i,2i+1$, and $2k+1$, such that the color $2k+1$ is only used for edges between vertices in $A_i$. Do so for each $i$. Now, the sets of colors we used for each $H_i$ were disjoint, save for the color $2k+1$. This was only used to color edges between points in $A_i$, though, so since the $A_i$'s are pairwise disjoint, this will not cause any conflicts. Thus in the end, we have a Borel edge coloring using the colors $1,2,\ldots,2k+1$, as desired.

\end{proof}

We now return to our proof of Theorem \ref{th:main} in the case $\Gamma/\Delta \cong \Z$. Fix a Borel $2k+1$-edge coloring $c$ of $H$, say using the colors $1,2,\ldots,2k+1$

Let $\gamma \in S_n$ for some $n > 0$. For each $\Delta$-orbit $y$, $\gamma$ corresponds to an edge from $y$ to $n \cdot y$ in $H$. Suppose $c$ assigns the color $l$ to that edge. Then let us give the edges $(x,\gamma \cdot x)$ in $G$ the color $l$ for each $x \in y$. Of course, since $x \neq x' \Rightarrow \gamma \cdot x \neq \gamma \cdot x'$, this does not cause any conflicts. Also, since $c$ was a coloring of $H$, doing this for all $\gamma \in S \setminus S_0$ does not cause any conflicts. 

It remains to color the edges corresponding to generators in $S_0$. These are the edges within each $\Delta$-orbit. For every such orbit $y$, the induced subgraph $G \res y$ is $d_0$-regular, so by Vizing's theorem it can be $d_0+1$-edge colored, say with the colors $2k+2,\ldots,2k+d_0+2 = d+2$. Since there are only finitely many such colorings for each orbit, we may choose one of them for each orbit in a Borel fashion. We have now $d+2$-edge colored our graph $G$ in a Borel fashion.

Finally, for each $\Delta$-orbit $y$, there must be some color $l \in \{1,\ldots,2k+1\}$ which does not appear on any edges incident to $y$ in our coloring $c$ of $H$, since $H$ is $2k$-regular. This means that, in our $d+2$-edge coloring above, none of the edges incident to a point in $y$ have the color $l$. Therefore, in the $d_0+1$-coloring of $y$ we have, we may replace the color $d+2$ with the color $l$ without causing any new conflicts. Doing so, we improve our coloring to a $d+1$-edge coloring, and so are done.

Finally, let us address the case $\Gamma/\Delta \cong D_\infty$. The argument which showed the sufficiency of Lemma \ref{lem:Z} was completely general, so here it suffices to show

\begin{lemma}
Theorem \ref{th:main} holds when $\Gamma = D_\infty$, even when multiplicity for generators is allowed. 
\end{lemma}

This can be proved similarly to Lemma \ref{lem:Z}. If $\gamma \in D_\infty$ is an order two element, no two of the edges it corrresponds to share an edge, So they can be Borel colored with a single color. Else, $\gamma$ has infinite order, so the edges corresponding to $\gamma$ and $\gamma\inv$ can be Borel 3-colored. As in the proof of Lemma \ref{lem:Z}, the third color here can be the same for every such $\gamma$, and used sparsely enough for each $\gamma$ so that there is no conflict in the end.

Thus, Theorem \ref{th:main} is proved. 

\end{proof}

\end{document}